\documentclass[preprint,12pt]{elsarticle}

\usepackage{color}
\usepackage{ulem}

\usepackage{amssymb}

\def\inf{{\rm{inf}}}
\def\mat{{\rm{mat}}}

\def\Vol{{\rm{vol}}}
\def\total{{\mathrm{total}}}
\def\cB{{\mathcal{B}}}

\def\Rconf{{\mathfrak{R}}}
\def\talpha{{t}}

\journal{Physica A}

\begin{document}

\begin{frontmatter}



\title{Fractal Structure of Equipotential Curves on a Continuum Percolation
Model}


\author[canon]{Shigeki Matsutani}
\ead{matsutani.shigeki@canon.co.jp}
\author[canon]{Yoshiyuki Shimosako}
\author[canon]{Yunhong Wang\corref{present}}

\address[canon]{Analysis technology development center,
Canon Inc. 3-30-2, Shimomaruko,
 Ohta-ku, Tokyo 146-8501, Japan}


\cortext[present]{[Present address]Microcreate System Co.,
2-1-15 4C, Chuou, Yamato, Kanagawa 242-0021, Japan}

\begin{abstract}
We numerically
investigate the electric potential distribution over
a two-dimensional continuum percolation model between the electrodes.
The model consists of overlapped conductive particles on the background
with an infinitesimal conductivity.
Using the finite difference method,
we solve the generalized Laplace equation
and show that in the potential distribution,
there appear the {\it{quasi-equipotential
clusters}} which approximately and locally
have the same values like steps and stairs.
Since the quasi-equipotential clusters have the fractal structure,
we compute the fractal dimension of equipotential curves
and its dependence on the volume fraction over $[0,1]$.
The fractal dimension in [1.00, 1.246] has a peak at
the percolation threshold $p_c$.

\end{abstract}


\end{frontmatter}



\section{Introduction}
In the series of articles \cite{MSW,MSW2},
we have studied the electric conductivity
of a percolation model using the finite difference method 
(FDM) \cite{LV}.
The purpose of these studies is to
reveal the electric properties 
over the continuum percolation model (CPM). 
As we explored three dimensional shape effects on the conductivity on 
CPMs in the previous articles \cite{MSW,MSW2}, 
in this article we investigate the fractal properties of the 
solutions of the two dimensional generalized Laplace equation,
\begin{equation}
\nabla\cdot \sigma \nabla \phi = 0,
\label{eq:0-1}
\end{equation}
for the conductivity distribution $\sigma(x,y)$ of 
CPM.
Recently materials consisting of the conductive particles 
in a base material with extremely high resistance are studied 
\cite{ASS, TMASC, KC}.
We handle the binary local conductivities $\sigma(x,y)$
of the binary materials in CPM, i.e., conductive particles
in the background with an infinitesimal conductivity.
The infinitesimal conductivity, i.e., the extremely high resistance,
avoids indeterminacy of the solutions of (\ref{eq:0-1}) all over there,
because approaching to zero differs from the zero itself.

Two dimensional CPMs
of the overlapped circles were studied well in Refs.\cite{Is} and
\cite{SA} and references therein.
In this article, we also deal only with CPM whose particles
are overlapped circles with the same radius.

As the solution of the generalized
Laplace equation (\ref{eq:0-1}) provides the electric
potential distribution, 
we show that
there exist {\it{quasi equipotential
clusters}}, in which the potential behaves like steps and
stairs approximately.
There are many pieces whose potentials are approximately flat because
their electric connections are very small and thus the resistances are so
high at the boundaries whereas the inner resistance of the pieces or the 
composite particles is very small.

The existence of the quasi-equipotential clusters means the
anomalous behaviors of the equipotential curves.
Though the relation between the percolation cluster and fractal
structure has been studied \cite{SA},  
in this article, we provide the result
how the fractal dimension of the equipotential curves
behaves depending on the volume fraction $p \in [0,1]$.
Under the percolation threshold $p_c$, our computation can read the
computation of dielectric behavior on a
random configuration of metal particles in
the dielectric matter, which was reported in Refs.~\cite{GB}
and \cite{SMKK} if we interpret the conductivity as
the dielectric constant.

Further recently, some classes of statistical models including several
two dimensional percolation models and the self-avoiding random walk
have investigated in the framework of 
Schramm-Loewner evaluation (SLE) and conformal field theory
\cite{C3,Law,Smir,LSW1,LSW2}.
Using their conformal properties, it is proved rigorously that 
some of the universal quantities such as the critical exponent and
the fractal dimension are expressed by rational numbers such as 5/4 and
4/3.
Since by taking a certain limit, the potential $\phi$ in (\ref{eq:0-1})
can show the conformal property, it is expected that
our numerical computational results should be also interpreted
in their framework in future. We give some comments in
Sec.3.3.

Contents in this article are as follows.
Sec.2 gives our computational method;
Sec.~2.1 is on the geometrical setting in CPMs and
Sec.~2.2 gives the computational method of the conductivity
over CPMs using FDM, which 
is basically the same as that in the previous article \cite{MSW}.
Sec.~3 shows our computational results of the solutions of 
(\ref{eq:0-1}) in CPMs, and we discuss
 our results from physical viewpoints there.

\bigskip

\section{Computational method}

In this section, we explain our computational method of CPM. 
Though we deal only with two dimensional case,
its essential is the same as the three dimensional case
 using FDM \cite{LV} whose detail is in  the previous article \cite{MSW}.


\subsection{Geometrical setting}

We set particles parametrized by their positions $(x,y)$ 
into a box-region $\cB:=[0,x_0]\times [0,y_0]$ 
at random and get a configuration $\Rconf_n$ as one of 
CPMs. In this article, we set $x_0 = y_0 = L_a$, 
and handle three boxes $\cB_a$ $(a = 1, 2, 3)$ with different size,
$L_1 = 40.96$, $L_2 = 81.92$, and $L_3 = 163.84$.
The particle corresponds to a stuffed circle with
the same radius $\rho=1$,
$B_{x_i, y_i}:=
\{ (x, y) \in \cB \ | \ |(x,y)-(x_i, y_i)| \le \rho\}$.
The configuration $\Rconf_n$ is given by 
$\Rconf_n:=\bigcup_{i=1}^n B_{x_i, y_i}$,
where 
each center $(x_i,y_i)$ is given at uniform random in $\cB$.
We allow their overlapping.

By monitoring the total volume fraction which is
a function of $\Rconf_n$ and is denoted by $\Vol(\Rconf_n)$, 
we continue to put the particles as long as
$\Vol(\Rconf_n) \le p$ for the given volume fraction $p$.
We find the step $n(p)$ such that
$\Vol(\Rconf_{n(p)-1}) \le p$ and
$\Vol(\Rconf_{n(p)}) > p$.
The difference between $\Vol(\Rconf_{n(p)-1})$ and 
$\Vol(\Rconf_{n(p)})$ is at most 
$1.87\times 10^{-3}$ for $\cB_1$,
$4.68\times 10^{-4}$ for $\cB_2$, and
$1.17\times 10^{-4}$ for $\cB_3$,
we regard $\Vol(\Rconf_{n(p)})$ as each $p$ hereafter
under this accuracy.

\begin{figure}[htbp]
\begin{center}
\includegraphics[width=12cm]{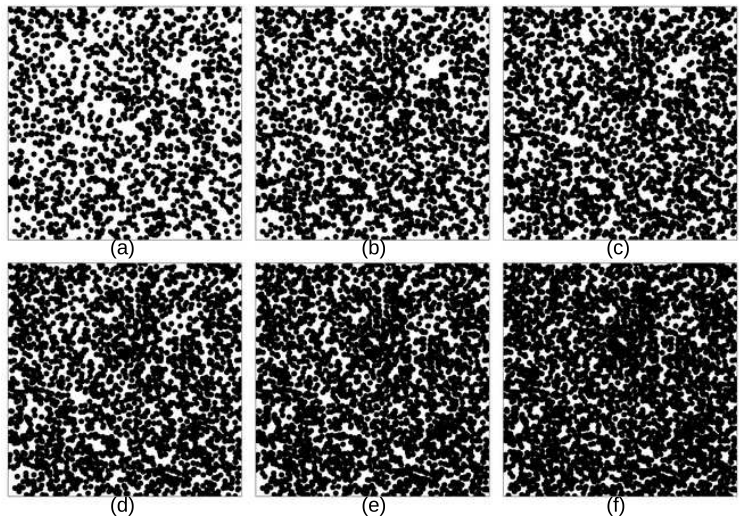}
\end{center}
\caption{The configurations of the seed $i_s=1$ with
the volume fractions 0.5, 0.6, 0.65, 0.7, 0.75, 0.8 for
(a), (b), (c), (d), (e), (f) respectively. }
\label{fig:sigma_config}
\end{figure}

Since we use the pseudo-randomness to simulate the random
configuration $\Rconf_{n(p)}$ for given $p$,
the configuration $\Rconf_{n(p)}$ depends upon the
 seed $i_s$ of the pseudo-random which we choose and thus
let it be denoted by $\Rconf_{p, i_s}$.
For the same seed $i_s$ of the pseudo-random,
a configuration $\Rconf_{p,i_s}$ of a volume fraction $p$ 
naturally contains a configuration $\Rconf_{p', i_s}$ of $p'<p$
due to our algorithm.
Hence the elements in the set of the configurations 
$\{\Rconf_{p,i_s}\ | \ p \in [0,1]\}$ with
the same seed $i_s$ are relevant.
Fig.\ref{fig:sigma_config} illustrates the configurations  
of the seed $i_s=1$ for several $p$.

\subsection{Computation of the potentials in CPM}

To apply FDM to the generalized Laplace equation in CPM,
we use three lattices, $1024 \times 1024$, 
$2048 \times 2048$, and $4096 \times 4096$, to represent the 
$\cB_1$, $\cB_2$ and $\cB_3$ respectively.
The radius of the particle $\rho=1$ corresponds to $25$ meshes.

Further we set the binary local conductivity $\sigma(x, y)$
which consists of the conductive particles 
with the conductivity density $\sigma_\mat = 1$,
and the background with
infinitesimal conductivity $\sigma_{\inf}=10^{-4}$.

In order to compute the potential distribution $\phi(x,y)$,
we set $\phi=1$ and $\phi=0$ on the upper and
the lower segments respectively as the boundary condition
 corresponding to the electrodes.
As the side boundary condition, we use the natural boundary
for each  $y$-boundary so that
current normal to each boundary vanish.

Following our algorithm of FDM as  mentioned in detail in Ref.\cite{MSW},
we numerically solve the generalized Laplace equations 
(\ref{eq:0-1})
for the conductivity distribution $\sigma(x,y)$ of each 
$\Rconf_{p, i_s}$.

Then we obtain the total conductivity
$\sigma_\total$
of the system after we integrate the current $\sigma \nabla \phi$
over the line parallel to the $x$-axis.
Since $\sigma_\total$ is determined for each $\Rconf_{p,i_s}$,
$\sigma_\total$ is a function of the volume fraction $p$
and the seed $i_s$ of the pseudo-random. We denote it by $\sigma_\total(p)$.

The dependence of the total conductivities on
the volume fraction as a conductivity curve 
is given by
\begin{equation}
\sigma_\total(p) = 
\frac{(p - p_c)^\talpha}{(1 - p_c)^\talpha}
\theta(p-p_c),
\label{eq:sigmatotal}
\end{equation}
where $p_c$ is the threshold, $\talpha$ is the 
critical exponent, or merely called exponent,
and $\theta$ is a Heaviside function, i.e.,
$\theta(x) = 1$ if $x \ge 0$ and vanishes otherwise.
(\ref{eq:sigmatotal}) is  
represented by a difference between the holomorphic
and anti-holomorphic functions in terms of
 the Sato hyperfunction theory \cite{Im}.
As mentioned in Introduction,
the conformal properties, which are given by
the holomorphic and anti-holomorphic functions, 
in a two-dimensional percolation model 
are recently revealed \cite{C3,L,Smir}, whereas
the holomorphic function, in general, determines its global behavior
from its local behavior because of the identity theorem
in theory of the complex analysis.
By considering these facts, we assume, in this article, that 
(\ref{eq:sigmatotal}) is defined over $[0,1]$ 
rather than the restricted region around $p_c$.
At least, from a viewpoint of numerical computations,
the least square error method shows 
the good descriptions of the conductivity curves 
in terms of the equation (\ref{eq:sigmatotal}) over $[0, 1]$.
Since the threshold  $p_c$ and the exponent
$\talpha$ are determined by
$\{\Rconf_{p,i_s} \ | \ p \in [0,1]\}$, they also 
depend on the seed $i_s$ of the pseudo-random.

\bigskip

\section{The results and discussions}

\subsection{Conductivity curve}

The ordinary computations of the conductivity
over the percolation system are performed using
the Y-$\Delta$ algorithm \cite{JJHS} by assuming that 
the conductivity of the background insulator vanishes \cite{GTDD};
(\ref{eq:0-1}) is not defined over there 
in the ordinary approaches.

In our computations using FDM,
we handle the binary conductivity distribution $\sigma(x,y)$ 
due to the demands of the recent material 
science \cite{ASS,TMASC,KC}.

\begin{figure}[htbp]
\begin{minipage}{0.5\hsize}
 \begin{center}
 \includegraphics[height=7cm, angle=270]{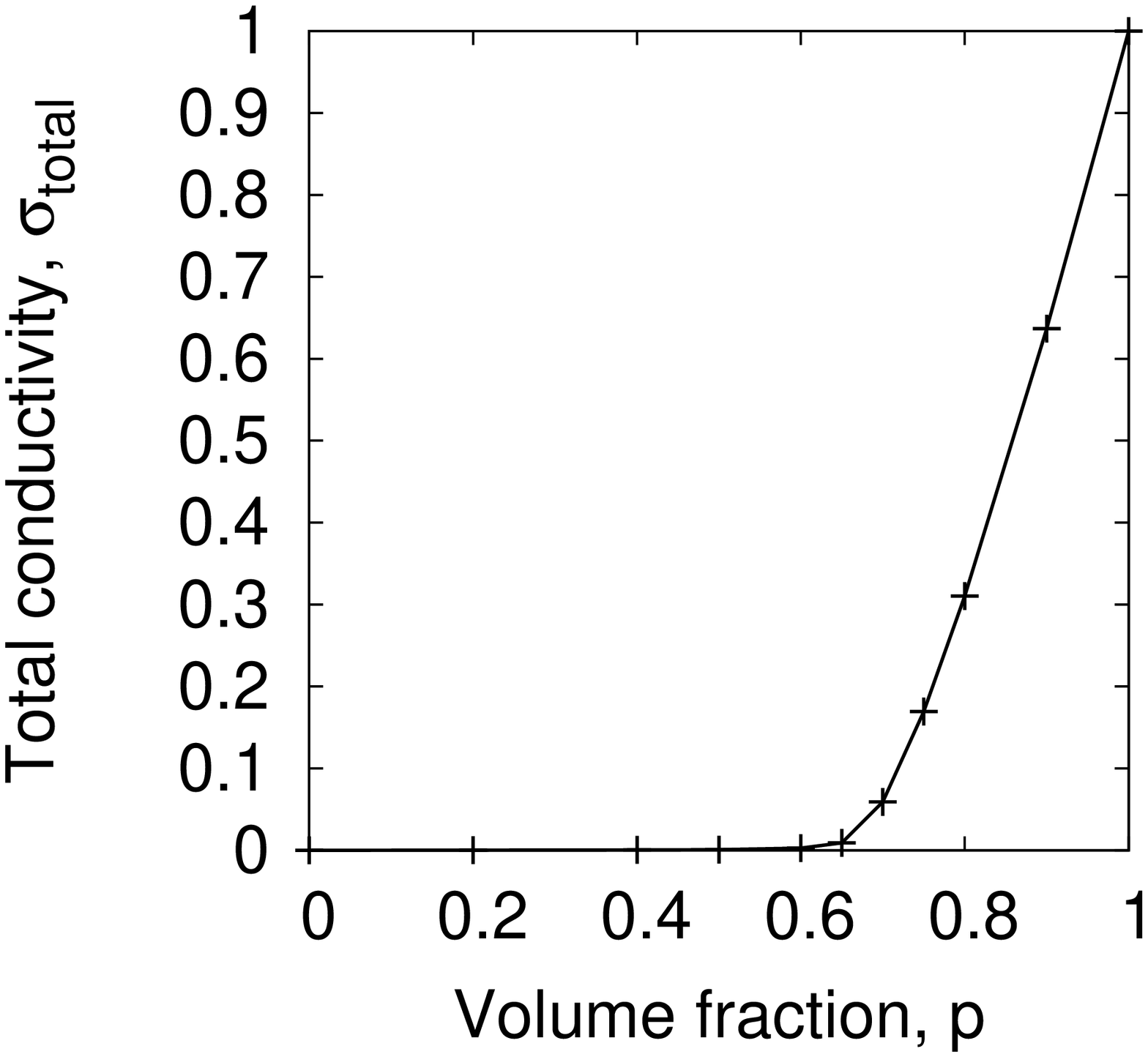} \newline
 (a)
 \end{center}
\end{minipage}
\begin{minipage}{0.5\hsize}
 \begin{center}
 \includegraphics[height=7cm, angle=270]{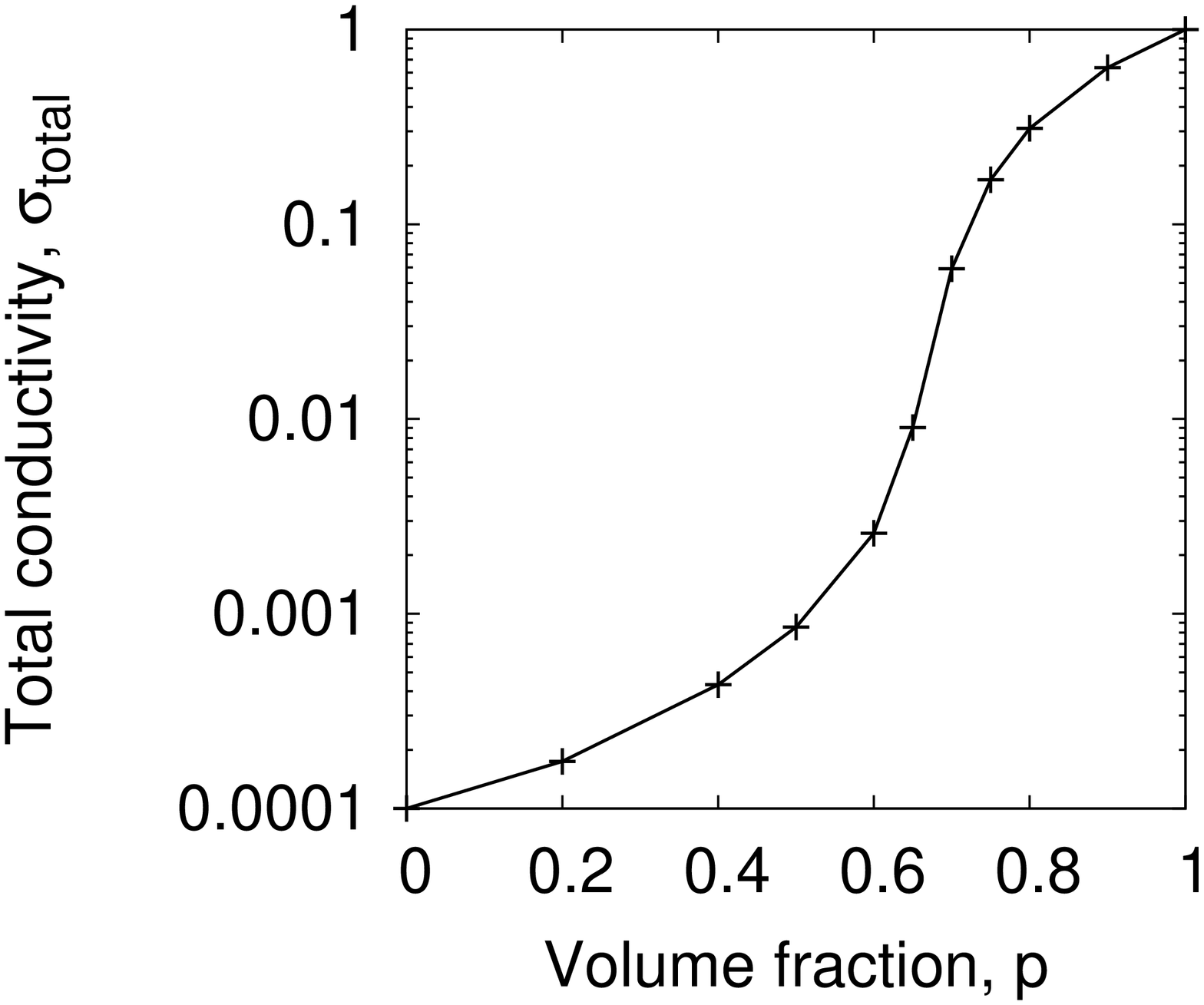} \newline
 (b)
 \end{center}
\end{minipage}
\caption{
Average of
conductivity curves for $\cB_3$;  
(a) for the linear total conductivity $\sigma_\total$ 
and (b) for its logarithm scale.
}
\label{fig:conductivity_curve} 
\end{figure}

Even though we deal with the binary conductivity distribution 
with finitely different resistances, 
the conductivity curve is naturally obtained as in 
Fig.\ref{fig:conductivity_curve}.

Fig.\ref{fig:conductivity_curve}(a) exhibits the linear scale 
behavior of 
average of the total conductivities 
in $\cB_3$
whereas Fig.\ref{fig:conductivity_curve}(b) shows its
logarithm property. Even though Fig.\ref{fig:conductivity_curve}(b) 
illustrates the property of the binary materials, the linear scale
behavior of the total conductivity is described well 
by (\ref{eq:sigmatotal})
because the ratio between the
binary conductivities $\sigma_\mat$ and $\sigma_\inf$ 
is sufficiently large.

We computed six cases with different seeds $i_s$ of the 
pseudo-randomness for the three boxes, $\cB_1$, $\cB_2$, and 
$\cB_3$ and determine the thresholds $p_c$ and 
the exponents $\talpha$ using the least mean square method
respectively.
\begin{figure}[htbp]
\begin{minipage}{0.5\hsize}
 \begin{center}
 \includegraphics[height=8cm, angle=270]{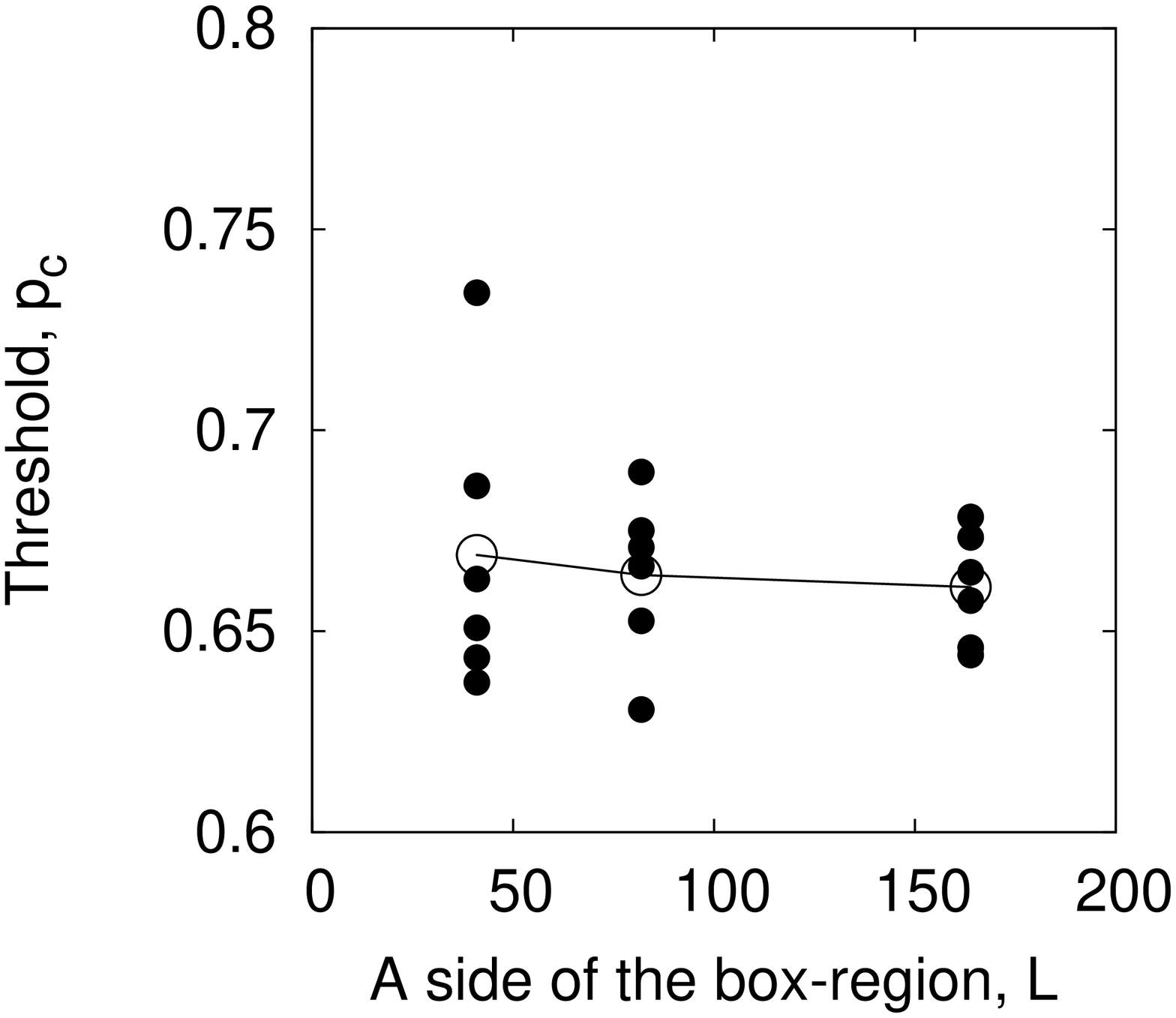} \newline
 (a)
 \end{center}
\end{minipage}
\begin{minipage}{0.5\hsize}
 \begin{center}
 \includegraphics[height=8cm, angle=270]{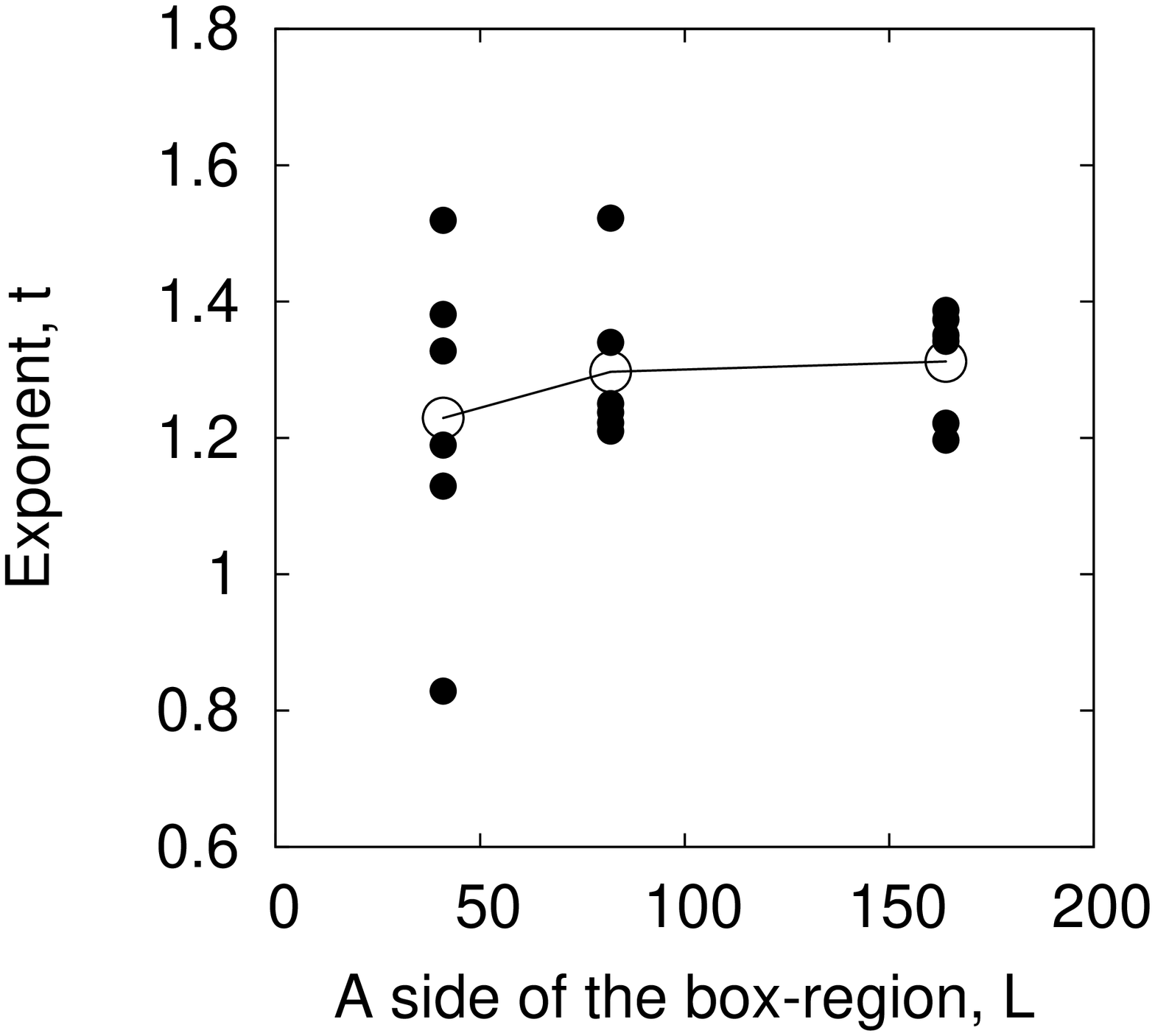} \newline
 (b)
 \end{center}
\end{minipage}
\caption{Threshold and Exponent vs. Box-size:
The filled circles show the individual values whereas
the open circles correspond to their averages.
}
\label{fig:lattice-dep}
\end{figure}
The obtained results are illustrated in Fig.\ref{fig:lattice-dep} and
Table \ref{table:lattice-dep}. Fig.\ref{fig:lattice-dep} shows
the dependence of the threshold and the exponent on the size of boxes.
The larger size of boxes is, 
the smaller the fluctuations of the threshold and the exponent are.
These results are also listed 
in Table.\ref{table:lattice-dep}.
\begin{table}[htbp]
\caption{The Size dependence of the threshold and the exponent:}
\label{table:lattice-dep}
\begin{center}
\begin{tabular}{|c|c|c|c|c|c|c|c|c|c|c|}
\hline
\multicolumn{1}{|c|}{ } &
\multicolumn{2}{|c|}{Threshold}&
\multicolumn{2}{|c|}{Exponent} \\
\hline
        &Average & Max-Min & Average & Max-Min  \\
\hline
$\cB_1$ & 0.669  & 0.097   & 1.229   & 0.690 \\
$\cB_2$ & 0.664  & 0.059   & 1.297   & 0.312 \\
$\cB_3$ & 0.661  & 0.032   & 1.312   & 0.190 \\
\hline
\end{tabular}
\end{center}
\end{table}


The average of the percolation threshold is 
$0.66\pm 0.02$ for $\cB_3$ 
which is equal to  $0.6763475(6)$ 
reported in Ref.~\cite{QZ} within our accuracy and
the exponent is given as $1.31\pm 0.10$ for $\cB_3$ 
which approximately recovers
the universal value $4/3$ of two-dimensional percolation \cite{Is}
within this accuracy.
Since the effect of the finite size of the analyzed region
makes the threshold smaller, 
the difference comes from the finite size effect \cite{MSW,MSW2}.
The direct computation of the total conductivity using FDM
provides the potential distribution as shown in 
Figs.\ref{fig:potentialdist} and \ref{fig:contours},
but it is for a finite region. Though it is very important that 
we evaluate the magnitude of the fluctuation, it is difficult
to remove these fluctuations in Fig.\ref{fig:lattice-dep}
in our FDM computations.
However in this article, we investigate the properties
of the potential distributions and these fluctuations do not 
essentially affect our results which are shown as follows.

\bigskip

\subsection{Potential distribution}

\begin{figure}[htbp]
\begin{center}
\includegraphics[width=12cm]{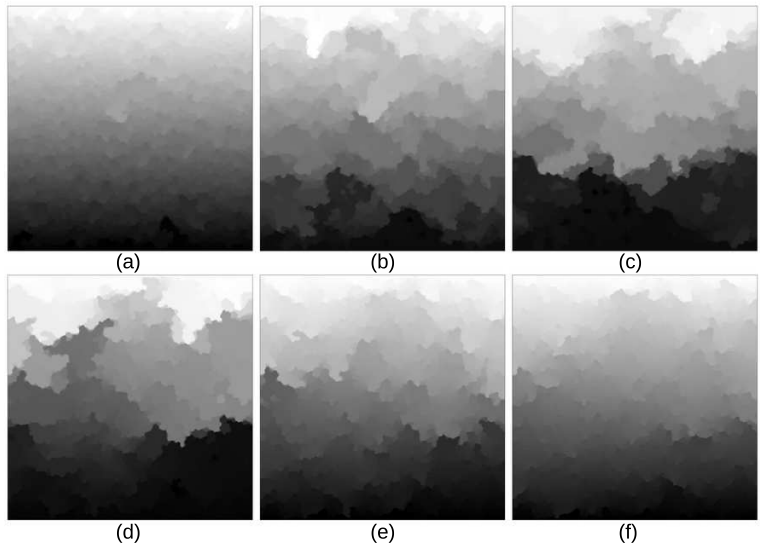}
\end{center}
\caption{Potential distributions of the seed $i_s = 1$
with the volume fractions 0.5, 0.6, 0.65, 0.7, 0.75, 0.8 for
(a), (b), (c), (d), (e), (f) respectively.
The black corresponds to $\phi=0$ whereas the white does
to $\phi=1$. The graduated gray interpolates them.
The quasi-potential clusters exist there.}
\label{fig:potentialdist}
\end{figure}

\begin{figure}[htbp]
\begin{center}
\includegraphics[width=12cm]{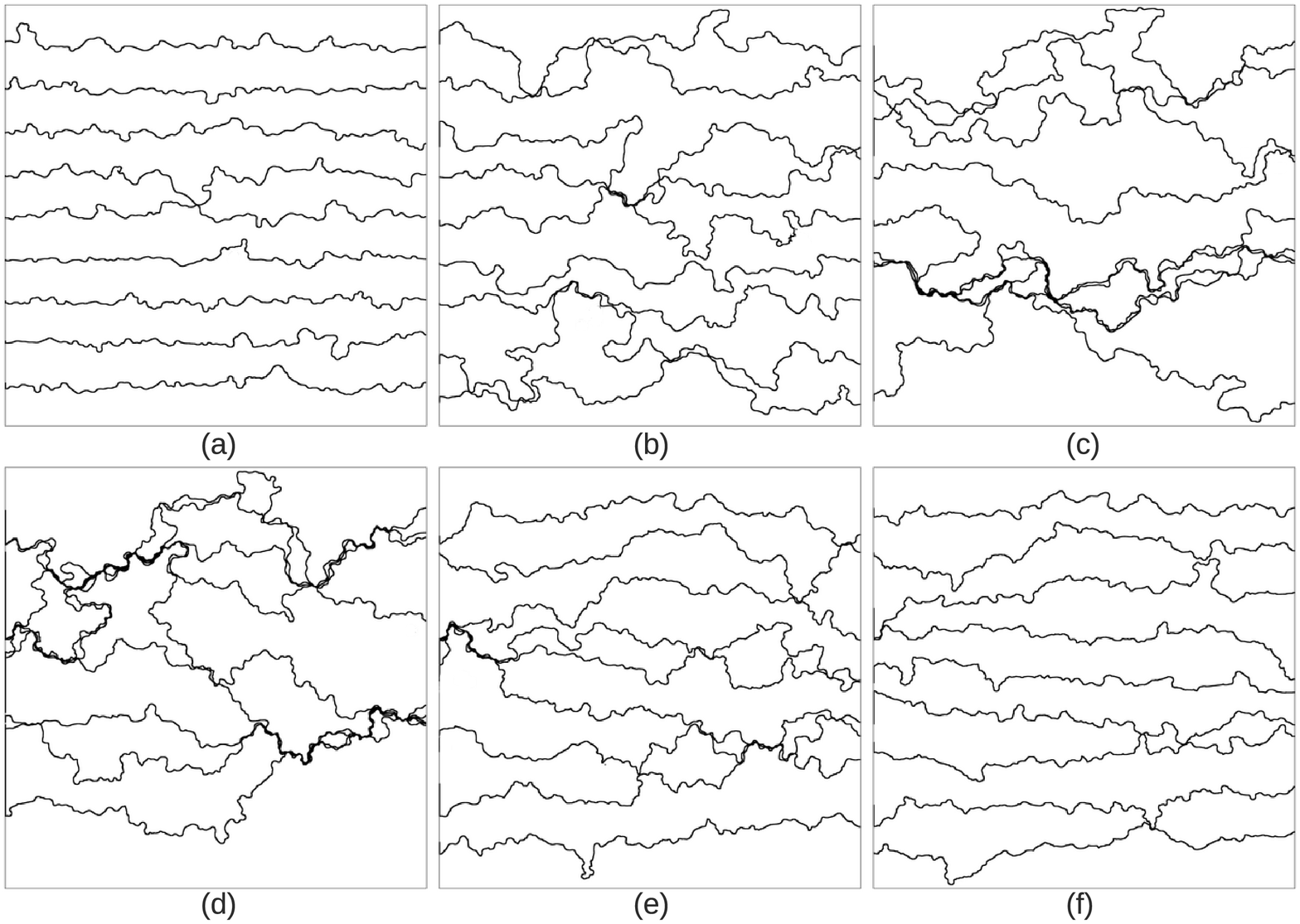}
\end{center}
\caption{Equipotential curves of the potential distributions
of the seed $i_s = 1$ with the volume fractions 0.5, 0.6, 0.65, 
0.7, 0.75, 0.8 for (a), (b), (c), (d), (e), (f) respectively. 
The curves correspond to the values $\phi = $ 0.1, 0.2, 
$\ldots$, 0.8, and 0.9.  The interval of the curves $\delta \phi$ is $0.1$.
}
\label{fig:contours}
\end{figure}

By numerically solving (\ref{eq:0-1}), we obtain the 
potential distribution $\phi(x,y)$ for each volume
fraction $p \in [0,1]$ and the seed $i_s$ of the pseudo-randomness 
as we display the results in Fig.\ref{fig:potentialdist} for $i_s = 1$.
Fig.\ref{fig:contours} shows its equipotential curves for 
$\phi = 0.1, 0.2, \ldots, 0.9$ of $\cB_3$ with $i_s =1$.
The interval of the equipotential curves $\delta \phi$ equals $0.1$.


\begin{figure}[htbp]
\begin{minipage}{0.5\hsize}
 \begin{center}
  \includegraphics[height=4cm]{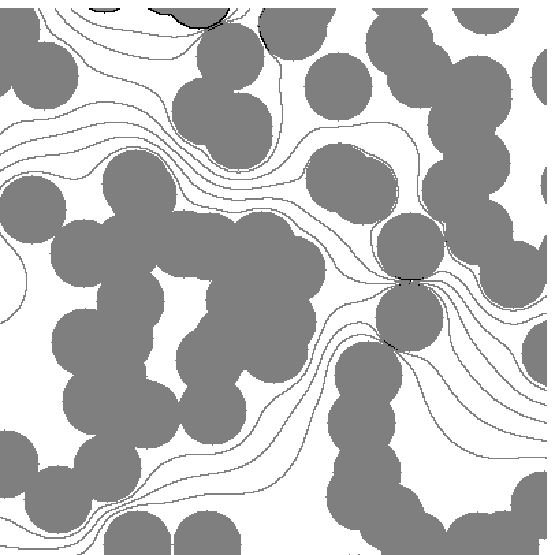} \newline
  (a)
 \end{center}
\end{minipage}
\begin{minipage}{0.5\hsize}
 \begin{center}
  \includegraphics[height=4cm]{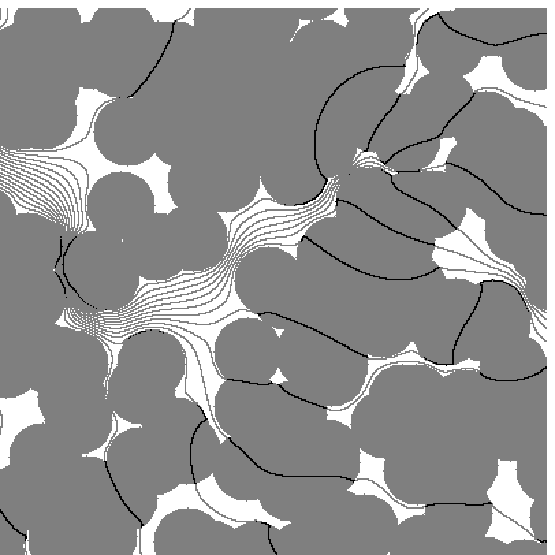} \newline
  (b)
 \end{center}
\end{minipage}
\caption{The equipotential curves with $\delta \phi = 0.0125$, 
and the configurations of particles;
(a) for $p=0.5$ and (b) for $p=0.8$.}
\label{fig:zoomin}
\end{figure}

Fig.\ref{fig:zoomin} is parts of Fig.\ref{fig:contours} with the 
configurations of
particles at $p = 0.5 < p_c$ and $p = 0.8 > p_c$.
The contours basically 
run to avoid the conductive particles due to the difference
of the local conductivity between $\sigma_\mat$ and $\sigma_\inf$.
Especially in Fig.\ref{fig:zoomin}(a) for 
the case $p=0.5$, the avoidance explicitly appears.

On the other hand, at $p=1$, the potential is simply described by 
$\phi(x,y) = y/y_0$ and the contours penetrate into the
conductive particles. Thus even for $p = 0.8$, the contours
penetrate into the conductive particles.
For the percolation clusters which are connected with both anode and cathode,
the contours cross the clusters.
On the other hand, around the approximately isolated 
percolation cluster which is 
not connected with electrodes approximately, the contours go through the
(narrow) gaps among the conductive
particles. Fig.\ref{fig:zoomin}(b) contains 
the  situations, in which
the contours partially penetrate into the conductive particles 
and partially avoid the conductive particles. However due to 
the above considerations,
these behaviors in Fig.\ref{fig:zoomin}(b) can be naturally interpreted.

Further Fig.\ref{fig:zoomin} shows that
the local equipotential curves depend upon 
the local configurations of these particles. For  given resolution
and local region,
the sets of particles surrounding the region locally play the role of
 the electrodes and the direction of the locally averaged gradient, i.e.,
the local electric field, strongly depends upon the resolutions and
upon the region which we average. 
For the locally averaged electric field within
 an appropriate scale $(L \gg a \gg 1)$,
we may represent its local direction by O(2) value spin, 
and thus our model could be
partially interpreted by XY-spin model with the external field
coming from the electrodes.
Since XY-model which is the O(2) valued spin system has the
fractal structure \cite{KO}, 
the randomness also brings the self-similar structure into our system.

Under the threshold $p_c$,
this phenomenon is essentially the same as the one 
of dielectric properties 
related to the electric breakdown, which
were studied by Gryure and Beale in Ref.~\cite{GB}
and recently by Stoyanov, Mc Carthy, Kollosche and Kofod 
in Ref.~\cite{SMKK},
though the authors handled the hard core model in which
the overlapping is prohibited instead of our soft core CPMs.
They studied the dielectric materials as background
with infinite dielectric particles (metal particles), which
are governed by the same generalized Laplace equation (\ref{eq:0-1}), 
though it is for the region of $p\in[0,p_c]$.
Fig.\ref{fig:zoomin}(a) corresponds to the
Fig.\ref{fig:sigma_config} in Ref.~\cite{GB} if the conductivity reads the
dielectric constant. 

In our computational condition on the conductivity,
we can determine the potential distributions over $\cB$ for 
 every volume fraction $p \in [0,1]$ due to the infinitesimal
conductivity $\sigma_\inf$.

As objects in the percolation theory exhibit 
the fractal structures related to the percolation clusters \cite{SA},
and to the breakdown \cite{PSIMMV,PD,FBS,RL},
the fractal dimension characterizes these properties.
We also computed the fractal dimensions 
but especially ones of the {\it{equipotential curves}}
in the following section.

\subsection{Fractal dimension of equipotential curves}

Fig.\ref{fig:potentialdist} provides the 256 degree 
shading pictures,
in which we find more complicate structures than those in
Fig.\ref{fig:contours}.
Fig.\ref{fig:potentialdist} illustrates that 
the generalized Laplace equation over a random configuration 
generates the self-similar structures.
Fig.\ref{fig:potentialdist} looks like pictures of the piled up 
mountains in distance using the monochrome water-painting.
There are many pieces whose potentials are approximately flat,
like leaves of the piled leaves.
Mandelbrot stated in Ref.~\cite{MB2} that
the fractal structures appear in many phenomena
including geographical geometry and these behaviors exhibit 
universal properties, which are profound mathematically 
beyond a shabby resemblance.
In fact, the recent studies in Refs.~\cite{C3,BB,Law}
show that the fractal structures are connected with 
wide fields in theoretical physics.
In fact, 
these properties are interpreted that there exist quasi-equipotential
clusters whose individual potential is approximately constant
and behaves like steps and stairs.
In other words,
the electric connections among the pieces
are very small and thus the resistance is so
high at the connections whereas the inner resistance of the piece or the 
composite particles is quite small.

Let us consider the effect of the quasi-equipotential clusters
on the equipotential curves in Figs.\ref{fig:contours} and \ref{fig:zoomin},
and the fractal dimensions 
of the {\it{equipotential curves}}.
On the computation of the fractal dimensions,
we used the box-counting method \cite{MB2}.
Since the curves in Fig.\ref{fig:contours} of the seed $i_s=1$
have the fractal structure, their fractal dimensions are
not trivial.
\begin{figure}[htbp]
\begin{minipage}{0.5\hsize}
 \begin{center}
  \includegraphics[height=8cm, angle=270]{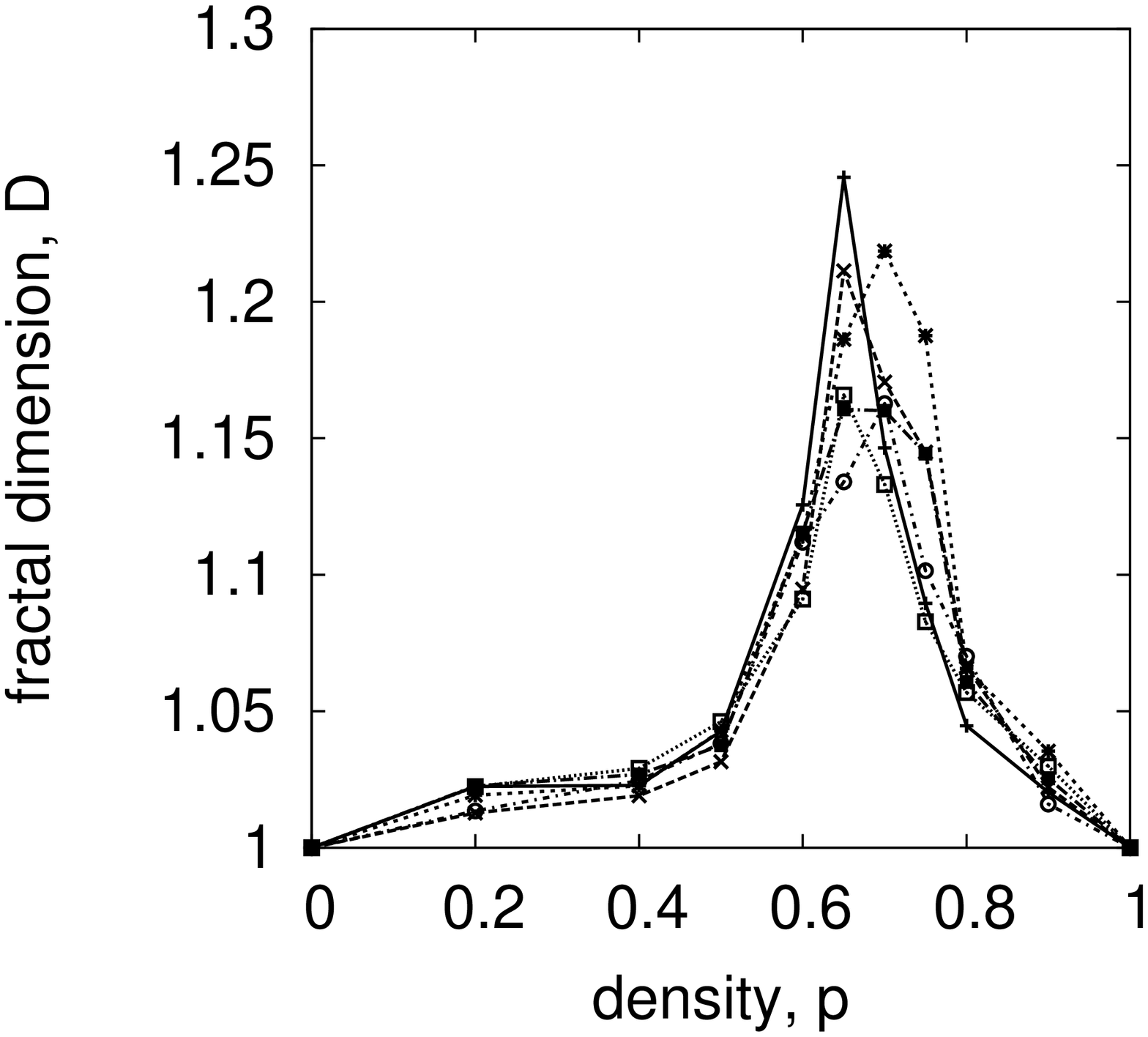} \newline
  (a)
  \includegraphics[height=8cm, angle=270]{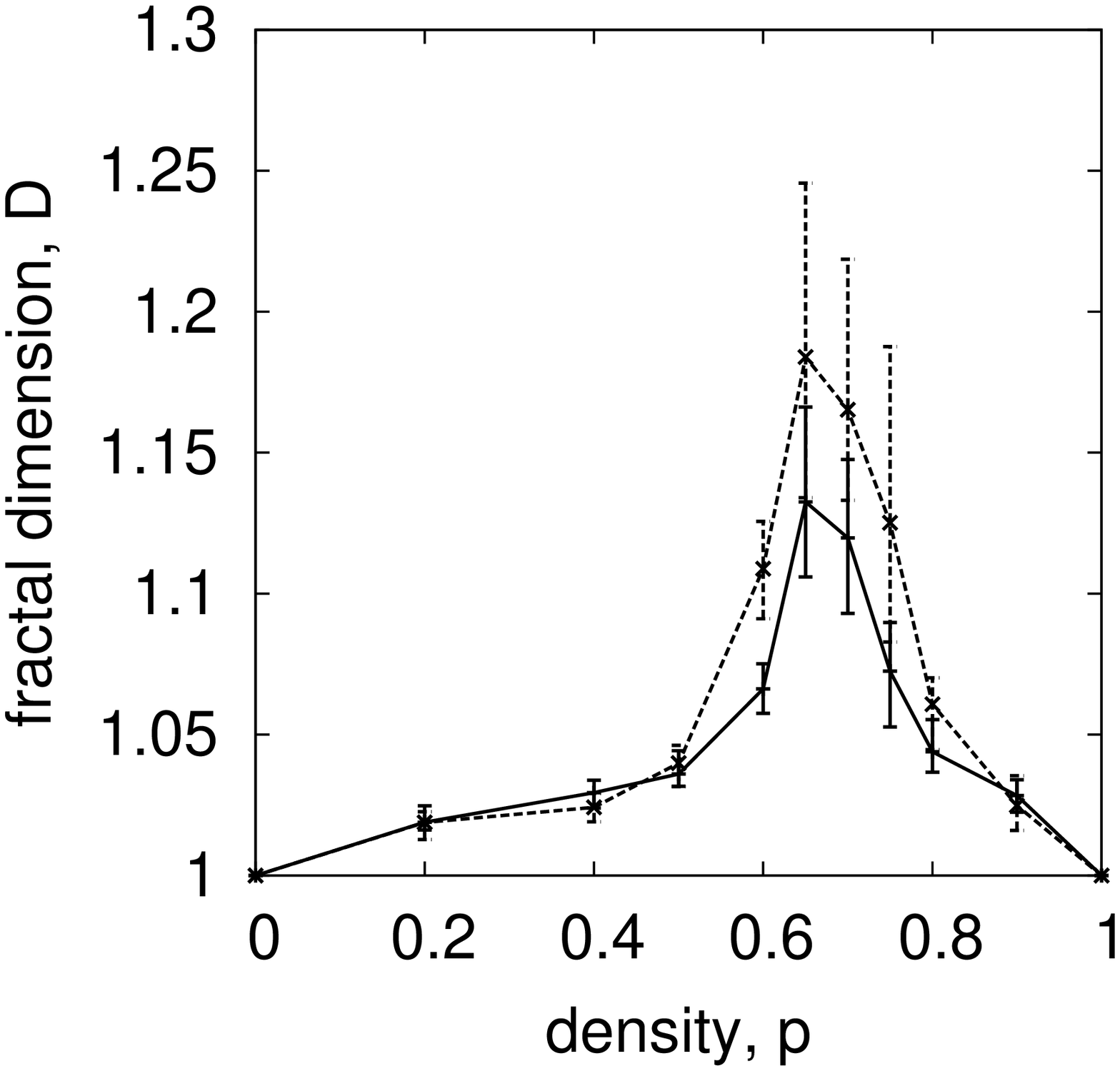} \newline
  (c)
 \end{center}
\end{minipage}
\begin{minipage}{0.5\hsize}
 \begin{center}
  \includegraphics[height=8cm, angle=270]{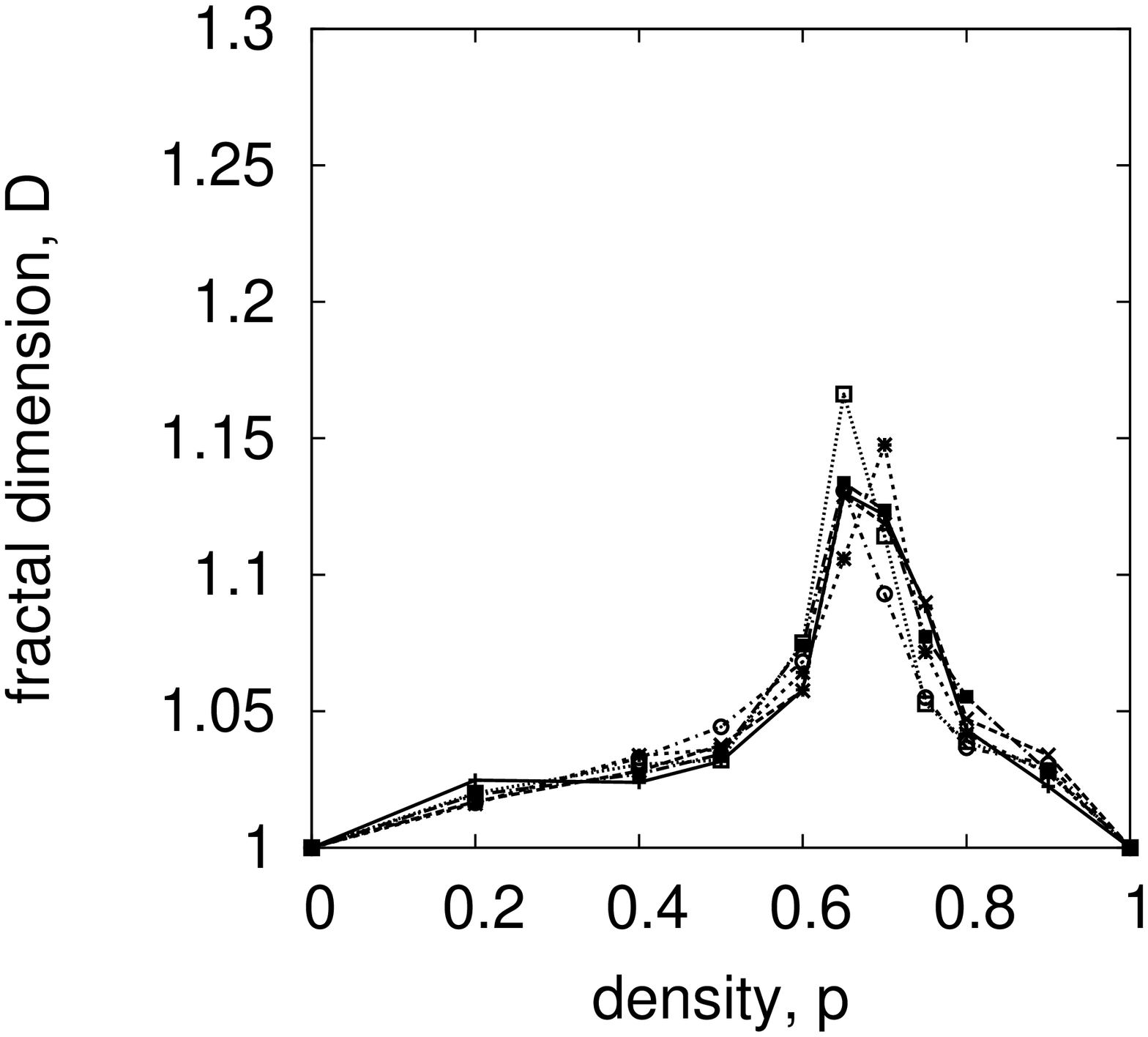} \newline
  (b)
  \phantom{\includegraphics[height=8cm, angle=270]{fig7c.eps} \newline}
  \phantom{(d)}
 \end{center}
\end{minipage}
\caption{The fractal dimension vs. volume fraction $p$ for six seeds:
(a) for $\cB_2$, (b) for $\cB_3$, 
and (c) for their average behaviors
with error bars. A dotted line and a solid line in (c) correspond to
$\cB_2$ and $\cB_3$ respectively.
}
\label{fig:fractaldim}
\end{figure}

The dependence of the fractal dimensions $D(p)$
of the equipotential curves 
with the interval of the equipotential curves $\delta \phi = 0.1$
on the volume fraction $p \in [0,1]$ for six different seeds $i_s$ is 
presented in Fig.\ref{fig:fractaldim};
Fig.\ref{fig:fractaldim} (a) shows the behavior of $\cB_2$, and
(b) shows that of $\cB_3$.
Fig.\ref{fig:fractaldim} (c) displays their average with error bars 
for their maximum and minimum values.
Since the percolation threshold $p_c$ of $\cB_2$ has larger fluctuation
than one of $\cB_3$, the fractal dimension $D(p)$ of $\cB_2$
also has larger fluctuations than ones of $\cB_3$.
However both behaviors are essentially the same without
depending on the size of the system.

Mandelbrot showed that the geographical curves like the coastlines
sometimes show fractal behavior \cite{MB2} and have
the fractal dimensions 1.0$\sim$1.25.
The equipotential curve is obviously a kind of self-avoiding curve
which is given in the randomized configuration,
though it is not clear whether the curve belongs to one of 
the universality classes, such as the self-avoiding random walk (SAW), 
the brownian frontier (BF), and the loop-erased self-avoiding walk (LESAW).
It is known that these fractal dimensions of SAW, BF, and LESAW are 4/3, 4/3 
and 5/4 respectively \cite{LSW2,MSN,LSW1,BB}. 
Further it is also known that 
the front of percolation cluster has 4/3 dimension as
results of numerical computations \cite{Is,Sap}.

Our fractal dimension $D(p)$ is in [1.00, 1.246] and has a peak at
the threshold $p_c$.
In other words, around $p_c$, the equipotential curves behave complicatedly
and have the self-similar structures.
On the other hand,
for the $p = 0$ and $p=1$ cases, the dimension must be equal to $1$.
Due to larger fluctuation of the threshold,
the curves in Fig.\ref{fig:fractaldim} (a) looks rude.
However these properties that the curve has the endings at $1$ and
has a peak around $p_c$ are universal.
Further the crudeness causes the values at the peaks larger than ones
of Fig.\ref{fig:fractaldim} (b). 
The larger the size of region is, the sharper the peak is,
as shown in Fig.\ref{fig:fractaldim} (c).
It is difficult to find the peak value. Further
the computation of (\ref{eq:0-1}) at the 
threshold $p_c$ basically needs a more rigorous numerical treatment  than
others because the singular behavior decrease the accuracy
of the numerical computations. We carefully computed the value
(using fine-tuned convergence parameters, and $\sigma_{\inf} = 10^{-6}$)
in $\cB_3$ at $\overline{p_c}= 0.661$,
which is 1.255.

Further when we approach the threshold $p=p_c$ from below, 
the percolation cluster connecting with the electrode 
faces each other and they are nearly connected. 
Thus some of the equipotential curves overlap at the front of these 
percolation clusters as we partially find the overlaps 
with respect to the resolution in Fig \ref{fig:contours}.
Similarly isolated percolation clusters also contribute
the shape of the equipotential curves as the front of
quasi equipotential clusters.
In other words, the equipotential curves, at least partially,
 are expected to agree with the front of the percolation clusters whose
the fractal dimension is known as 4/3.
Thus it is expected that the fractal dimension $D(p)$ should have
a peak at the percolation threshold  $p = p_c$ whose value is larger
than 1.246 and 1.255.
Ziff conjectures that the dimension may be equal to $4/3$ at $p_c$ \cite{Ziff}.
Since in the limit $\sigma_\inf \to +0$, 
the equation (\ref{eq:0-1}) over $p> p_c$ is reduced to the
mixed Dirichlet-Neumann boundary problem of the 
Laplace
equation with a random boundary,
and $\phi$ is expressed by a real part of a holomorphic function,
it is expected that our computational results might be also interpreted
in the framework of SLE and conformal field theory.

\section*{Acknowledgment}
We thank Professor Robert M. Ziff for giving our attention to
Ref.\cite{QZ} and helpful comments,
and the anonymous referee for his/her valuable comments.





\bibliographystyle{model1-num-names}
\bibliography{MSW120524}







\end{document}